\documentclass[10pt]{amsart}
\usepackage{tipa}
\usepackage{amssymb}
\usepackage{stmaryrd}
\usepackage{amsmath,amsfonts,amsthm,mathrsfs,bbm,array,subfigure}

\newtheorem{theorem}{Theorem}[section]

\newtheorem{corollary}[theorem]{Corollary}

\theoremstyle{definition}

\newtheorem{question}[theorem]{Question}

\theoremstyle{remark}

\numberwithin{equation}{section}

\begin{document}

\title[On the Diophantine equations $z^2=f(x)^2 \pm f(y)^2$]
{On the Diophantine equations $z^2=f(x)^2 \pm f(y)^2$ involving
Laurent polynomials}

\author{Yong Zhang}
\address{School of Mathematics and Statistics, Changsha University of Science and Technology,
Changsha 410114, People's Republic of China}
 \email{zhangyongzju@163.com}

\thanks{This research was supported by the National Natural Science Foundation of China (Grant No.~11501052).}

\author{Arman Shamsi Zargar}
\address{Independent Researcher, Ardabil, Iran}
 \email{shzargar.arman@gmail.com}

\subjclass[2010]{Primary 11D72, 11D25; Secondary 11D41, 11G05}

\date{  }

\keywords{Diophantine equation, Laurent polynomial, rational
parametric solution}

\begin{abstract}
By the theory of elliptic curves, we study the nontrivial rational
parametric solutions and rational solutions of the Diophantine
equations $z^2=f(x)^2 \pm f(y)^2$ for some simple Laurent
polynomials $f$.
\end{abstract}

\maketitle

\section{Introduction}
In 2010, A. Togbe and M. Ulas \cite{Togbe-Ulas} considered the
rational solutions of the Diophantine equations
\begin{equation}\label{Eq11}
z^2=f(x)^2+f(y)^2
\end{equation}
and
\begin{equation}\label{Eq12}
z^2=f(x)^2-f(y)^2,
\end{equation}
for $f$ being quadratic and cubic polynomials. At the same year, B.
He, A. Togbe and M. Ulas \cite{He-Togbe-Ulas} further investigated
the integer solutions of Eqs. (\ref{Eq11}) and (\ref{Eq12}) for some
special polynomials $f$.

In 2016, Y. Zhang \cite{Zhang2018} studied the nontrivial rational
parametric solutions of the Diophantine equations
\[
f(x)f(y)=f(z)^n,
\]
for $n=1,2,$ concerning the Laurent polynomials $f(x)=ax+b+c/x$.

In this paper, we continue the study of \cite{Zhang2018} and
consider the nontrivial rational parametric solutions of Eqs.
(\ref{Eq11}) and (\ref{Eq12}) for some simple Laurent polynomials.
The nontrivial solution $(x,y,z)$ of Eqs. (\ref{Eq11}) and
(\ref{Eq12}) respectively means that $f(x)f(y)\neq0$ and $f(x)^2\neq
f(y)^2$.

Recall that a Laurent polynomial with coefficients in a field
$\mathbb{F}$ is an expression of the form
\[f(x)=\sum_{k}a_kx^k,a_k\in \mathbb{F},\]
where $x$ is a formal variable, the summation index $k$ is an
integer (not necessarily positive) and only finitely many
coefficients $a_k$ are nonzero. Here we mainly care about the simple
Laurent polynomials
\[f(x)=ax+b+\frac{c}{x},~ax^2+bx+c+\frac{d}{x},~ax+b+\frac{c}{x}+\frac{d}{x^2}\]
with nonzero integers $a,b,c,d$. In the sequel, without losing the
generality we may assume that $a=1$, guaranteed by the shape of Eqs.
(\ref{Eq11}) and (\ref{Eq12}).

By the theory of elliptic curves, we have

\begin{theorem} For $f(x)=x+b+c/x$ with nonzero integers $b,c$, Eq. (\ref{Eq11}) or Eq. (\ref{Eq12}) has infinitely many nontrivial rational parametric
solutions.
\end{theorem}

\begin{theorem} There are infinitely many Laurent polynomials of the form $f(x)=x^2+bx+c+d/x$ with nonzero integers $b,c,d$ such that
Eq. (\ref{Eq11}) or Eq. (\ref{Eq12}) has infinitely many nontrivial
rational solutions.
\end{theorem}

Noting that Corollary 1.3 is obtained via sending $x\mapsto 1/x$ in
Theorem 1.2.
\begin{corollary} There are infinitely many Laurent polynomials of the form $f(x)=x+b+c/x+d/x^2$ with nonzero integers $b,c,d$
such that Eq. (\ref{Eq11}) or Eq. (\ref{Eq12})  has infinitely many
nontrivial rational solutions.
\end{corollary}

\section{Proofs of the theorems}
\begin{proof}[\textbf{Proof of Theorem 1.1.}]
For $f(x)=x+b+c/x$, let \[x=T,y=tT.\] Then Eq. (\ref{Eq11}) equals
\[\begin{split}
T^2t^2z^2=&T^4t^4+2bT^3t^3+(T^4+2bT^3+2b^2T^2+4cT^2+2bcT+c^2)t^{2}\\
&+2bcTt+c^2=:g_1(t).\end{split}\]

Let $v=Ttz,$ and consider the curve $\mathcal{C}_1:~v^2=g_1(t)$. In
order to prove Theorem 1.1 we must show that the curve
$\mathcal{C}_1$ has infinitely many $Q(T)$-rational points. The
curve $C_1$ is a quartic curve with rational point $P=(0,c)$. By the
method of Fermat \cite[p. 639]{Dickson}, using the point $P$ we can
produce another point $P'=(t,v)$, which satisfies the condition
$tv\neq0$. Indeed, in order to construct a such point $P'$ we put
$v= pt^2 + qt +c$, where $p,q$ are indeterminate variables. Then we
have
\[v^2-g_1(t)=\sum_{i=1}^4A_it^i,\] where the quantities $A_i = A_i(p,q)$ are given by
\[\begin{split}
A_1=&2bcT-2cq,\\
A_2=&T^4+2bT^3+2b^2T^2+4cT^2+2bcT+c^2-2cp-q^2,\\
A_3=&2bT^3-2pq,\\
A_4=&T^4-p^2.
\end{split}\]
The system of equations $A_4=A_3=0$ in $p, q$ has a solution given
by\[p=-T^2,q=-bT.\]This implies that the equation
\[v^2-g_1(t)=\sum_{i=1}^4A_it^i=0\] has the rational roots $t=0$ and \[t=-\frac{4bcT}{T^4+2bT^3+(b^2+6c)T^2+2bcT+c^2}.\]
Then we have the point $P'=(t,v)$, where
\[\begin{split}
t=&-\frac{4bcT}{T^4+2bT^3+(b^2+6c)T^2+2bcT+c^2},\\
v=&c(T^8+4bT^7+(10b^2+12c)T^6+(12b^3+28bc)T^5\\
  &+(5b^4+28b^2c+38c^2)T^4+(12b^3c+28bc^2)T^3+(10b^2c^2+12c^3)T^2\\
  &+4bc^3T+c^4)/(T^4+2bT^3+(b^2+6c)T^2+2bcT+c^2)^2.
\end{split}\]

If $c\neq0,$ the discriminant of $g_1(t)$ is nonzero as an element
of $\mathbb{Q}(T)$, then $\mathcal{C}_1$ is smooth. By the method
described as in \cite[p. 77]{Mordell} (or \cite[p. 476, Proposition
7.2.1]{Cohen}), $\mathcal{C}_1$ is birationally equivalent to the
elliptic curve
\[\mathcal{E}_1:Y^2=X^3+B_1X+B_2,\]
where
\[\begin{split}
B_1=&27(T^8+4bT^7+(8b^2+8c)T^6+(8b^3+20bc)T^5\\
    &+(4b^4+12b^2c+30c^2)T^4+(8b^3c+20bc^2)T^3\\
    &+(8b^2c^2+8c^3)T^2+4bc^3T+c^4),\\
B_2=&54(T^4+2bT^3+(2b^2-2c)T^2+2bcT+c^2)(T^8+4bT^7\\
    &+(8b^2+14c)T^6+(8b^3+32bc)T^5+(4b^4+18b^2c+42c^2)T^4\\
    &+(8b^3c+32bc^2)T^3+(8b^2c^2+14c^3)T^2+4bc^3T+c^4).
\end{split}\]
Since the map $\varphi_1: \mathcal{C}_1\rightarrow \mathcal{E}_1$ is
complicated, we do not present the explicit equations for the
coordinates of it.

By the map $\varphi_1$, we get the point
\[P''=\varphi_1(P')=(X(P''),Y(P'')),\]
on the elliptic curve $\mathcal{E}_1$, in which
\[\begin{split}
X(P'')=&\frac{3}{4b^2T^2}(3T^8+12bT^7+(22b^2+36c)T^6+(20b^3+84bc)T^5\\
       &+(11b^4+52b^2c+114c^2)T^4+(20b^3c+84bc^2)T^3\\
       &+(22b^2c^2+36c^3)T^2+12bc^3T+3c^4),\\
Y(P'')=&-\frac{27}{8b^3T^3}(T^4+2bT^3+(b^2+6c)T^2+2bcT+c^2)(T^8+4bT^7\\
       &+(8b^2+12c)T^6+(8b^3+28bc)T^5+(3b^4+16b^2c+38c^2)T^4\\
       &+(8b^3c+28bc^2)T^3+(8b^2c^2+12c^3)T^2+4bc^3T+c^4).
\end{split}\]

A quick computation reveals that the remainder of the division of
the numerator of the $X$-coordinate of the point $P''$ by the
denominator $(2bT)^2$ is equal to
\[9c^3(4bT+c)\] with respect to $T$, which is thus non-zero provided $bc\neq0$. Then the $X$-coordinate of the point
$P''$ is not a polynomial. By a generalization of Nagell-Lutz
theorem (\cite[p. 268]{Connell}), for the elliptic curve over
$\mathbb{Q}(T)$ with the equation
$y^2=x^3+a(T)x^2+b(T)x+c(T)$,~$a(T),b(T),c(T)\in \mathbb{Z}[T]$,
points of finite order have coordinates in $\mathbb{Z}[T]$, then
$P''$ is of infinite order on $\mathcal{E}_1$. Hence, the group
$\mathcal{E}_1(\mathbb{Q}(T))$ is infinite.

Compute the points $[m]P''$ on $\mathcal{E}_1$ for $m=2,3,...$, next
calculate the corresponding point $\varphi_1^{-1}([m]P'') = (t_m,
v_m)$ on $\mathcal{C}_1$, then get various $\mathbb{Q}(T)$-rational
solutions $(x,y,z)$ of Eq. (\ref{Eq11}) for $f(x)=x+b+c/x$ with
nonzero integers $b,c$. For example, the point $P''$ on
$\mathcal{E}_1$ leads to the solution of Eq. (\ref{Eq11}):
\[\begin{split}
(x,y,z)=\bigg(&T,-\frac{4bcT^2}{T^4+2bT^3+(b^2+6c)T^2+2bcT+c^2},\\
              &(T^8+4bT^7+(10b^2+12c)T^6+(12b^3+28bc)T^5\\
              &+(5b^4+28b^2c+38c^2)T^4+(12b^3c+28bc^2)T^3\\
              &+(10b^2c^2+12c^3)T^2+4bc^3T+c^4)/(4bT^2(T^4\\
              &+2bT^3+(b^2+6c)T^2+2bcT+c^2))\bigg).
\end{split}\]
This completes the proof of Theorem 1.1 for Eq. (\ref{Eq11}). The
same method can be used to give a proof for Eq. (\ref{Eq12}).
\end{proof}

\begin{proof}[\textbf{Proof of Theorem 1.2.}]
For $f(x)=x^2+bx+c+d/x$, let \[y=tx.\] Then Eq. (\ref{Eq11}) becomes
\[\begin{split}
t^2x^2z^2=&(t^6+t^2)x^6+(2bt^5+2bt^2)x^5+(2ct^4+b^2t^4+2ct^2+b^2t^2)x^4\\
          &+(2dt^3+2bct^3+2dt^2+2bct^2)x^3+(4bdt^2+2c^2t^2)x^2\\
          &+(2cdt^2+2cdt)x+d^2t^2+d^2=:g_2(t),
\end{split}\]
The discriminant of $g_2(t)$ is $$-64d^2t^{10}(t-1)^6h(t)i(t),$$
where $i(t)$ is a polynomial of degree 12 in terms of $t$ and
\[\begin{split}
h(t)=&d^2t^6+(3d^2-bcd)t^5+(6d^2-5bcd+c^3+b^3d)t^4\\
     &+(7d^2-6bcd+2c^3+2b^3d-b^2c^2)t^3+(6d^2-5bcd\\
     &+c^3+b^3d)t^2+(3d^2-bcd)t+d^2.
\end{split}\]

To find rational $t$ such that Eq. (\ref{Eq11}) has infinitely many
rational solutions, we need $g(t)$ has a square factor of the form
$(\alpha x+\beta)^2$. Putting the discriminant of $h(t)$ be zero,
i.e.,
\[-d^4(d-bc)^2(27d^2-18bcd+4c^3+4b^3d-b^2c^2)^3(c^3-b^3d)^4=0,\]
and solving it for $d$, we obtain
\[d=bc,\frac{c^3}{b^3},\frac{9bc-2b^3\pm2\sqrt{-(3c-b^2)^3}}{27}.\]
Consider the case $d=bc,$ then
\[h(t)=c(t+1)^2(b^2t^2+c)(ct^2+b^2)\] The expression $h(t)$ vanishes whenever
\[t=-1,\pm\frac{\sqrt{-c}}{b},\pm\frac{b}{\sqrt{-c}}.\]
Consider the case $t=\frac{\sqrt{-c}}{b},$ and take $-c=k^2b^2,$
i.e., $c=-k^2b^2$, then $t=k$ and
\[f(x)=x^2+bx-k^2b^2-\frac{k^2b^3}{x}=\frac{(x+b)(x-bk)(x+bk)}{x}.\]
Then Eq. (\ref{Eq11}) reduces to
\[\begin{split}
x^2z^2=&(x+b)^2[(k^4+1)x^4+(-2bk^4+2bk^3)x^3+(b^2k^4-4b^2k^3-b^2k^2)x^2\\
&+(2b^3k^3-2b^3k^2)x+b^4k^2(k^2+1)].
\end{split}\]
To find infinitely many $x$ such that
\[\begin{split}
&(k^4+1)x^4+(-2bk^4+2bk^3)x^3+(b^2k^4-4b^2k^3-b^2k^2)x^2\\
&+(2b^3k^3-2b^3k^2)x+b^4k^2(k^2+1)
\end{split}\]
is a rational square, let $k^2+1=l^2,$ then
\[k=\frac{r^2-1}{2r},l=\frac{r^2+1}{2r},\] where $r\neq0,\pm1$ is an integer parameter. Then Eq.
(\ref{Eq11}) leads to
\[16r^4x^2z^2=(x+b)^2[a_4x^4+a_3x^3+a_2x^2+a_1x+a_0],\] where
\[\begin{split}
a_4&=r^8-4r^6+22r^4-4r^2+1,\\
a_3&=-2b(r^2-2r-1)(r-1)^3(r+1)^3,\\
a_2&=b^2(r^4-8r^3-6r^2+8r+1)(r-1)^2(r+1)^2,\\
a_1&=4b^3r(r^2-2r-1)(r-1)^2(r+1)^2,\\
a_0&=b^4(r-1)^2(r+1)^2(r^2+1)^2.
\end{split}\]
Hence, we need to study the rational points on the quartic curve
\[\mathcal{C}_2:~w^2=a_4x^4+a_3x^3+a_2x^2+a_1x+a_0=p(r).\] When $b\neq0$ and $r\neq 0,\pm1$, the discriminant of $p(r)$, i.e.,
\[\begin{split}
&16(17r^6-26r^5+27r^4+4r^3-5r^2-2r+1)(r^6+2r^5-5r^4-4r^3+27r^2\\
&+26r+17)(r^2+1)^4(r^2+2r-1)^4(r^2-2r-1)^4(r^2-1)^6b^{12}
\end{split}\]
is nonzero, then $\mathcal{C}_2$ is smooth. By the method described
as in \cite[p. 477]{Cohen}, and using Magma
\cite{Bosma-Cannon-Playoust}, one can observe that $\mathcal{C}_2$
is birationally equivalent to the elliptic curve
\[\begin{split}
\mathcal{E}_2:~Y^2=X^3+b_1X+b_2,
\end{split}\]
where
\[\begin{split}
b_1=&-b^4(13r^{12}+8r^{11}-70r^{10}+120r^9+371r^8-400r^7+140r^6+400r^5\\
&+371r^4-120r^3-70r^2-8r+13)(r^2-1)^2/3,\\
b_2=&2b^6(r^4+4r^3+18r^2-4r+1)(19r^{12}+8r^{11}-130r^{10}+120r^9+461r^8\\
&-400r^7+452r^6+400r^5+461r^4-120r^3-130r^2-8r+19)(r^2-1)^4/27.
\end{split}\]
Again, as the map $\varphi_2: \mathcal{C}_2\rightarrow
\mathcal{E}_2$ is complicated, we do not present the explicit
equations for the coordinates of it.

Let\[U=\frac{9X}{b^2},V=\frac{27Y}{b^3},\] then
\[\begin{split}
\mathcal{E}_{2}(r):~V^2=U^3+c_1U+c_2,
\end{split}\]
where
\[\begin{split}
c_1=&-27(13r^{12}+8r^{11}-70r^{10}+120r^9+371r^8-400r^7+140r^6+400r^5\\
&+371r^4-120r^3-70r^2-8r+13)(r^2-1)^2,\\
c_2=&54(r^4+4r^3+18r^2-4r+1)(19r^{12}+8r^{11}-130r^{10}+120r^9+461r^8\\
&-400r^7+452r^6+400r^5+461r^4-120r^3-130r^2-8r+19)(r^2-1)^4.
\end{split}\]

It is easy to see that $\mathcal{E}_{2}(r)$ contains the point
\[\begin{split}
Q=(&3(r^2-1)(r^6+4r^5+5r^4-8r^3+55r^2+4r-13),\\
&108(r^2+2r-1)(r^2-2r-1)(r^3-r+2)(r^2-1)^2).
\end{split}\]
By the group law, we have
\[\begin{split}
[2]Q=\bigg(&(3(3r^{16}-20r^{14}+16r^{13}+68r^{12}-64r^{11}-220r^{10}+400r^9+970r^8\\
           &-896r^7+1268r^6+688r^5+1428r^4-64r^3-516r^2-80r\\
           &+91))/(2(r^3-r+2))^2,\\
           &27(r^4+2r^2+5)(r^{12}+8r^9-3r^8+96r^6-48r^5+67r^4+64r^3+96r^2\\
           &-24r-1)(r^2+2r-1)^2(r^2-2r-1)^2/(2(r^3-r+2))^3\bigg).\end{split}\]

To prove that there are infinitely many rational points on
$\mathcal{E}_{2}(r)$, it is enough to find a point on
$\mathcal{E}_{2}(r)$ with nonintegral $U$-coordinate. When the
numerator of the $U$-coordinate of $[2]Q$ is divided by
$(2(r^3-r+2))^2$, the remainder equals
\[R=9(-724r^5+1529r^4-1876r^3-1206r^2+1896r-2287)\] and $R\neq 0$ when
$r\neq0,\pm1,$ so the $U$-coordinate of $[2]Q$ is not a polynomial.

For $r\neq 0,\pm1,3$ and $-1631\leq r\leq1626$ one can check that
$|R|/(2(r^3-r+2))^2$ is not an integer, and that it is nonzero and
less than 1 in modulus for $r>1626$ and $r<-1631$. Hence, for all
$r\neq0,\pm1,3$ the point $[2]Q$ has nonintegral $U$-coordinate and
hence, by the Nagell-Lutz Theorem (\cite[p. 56]{Silverman-Tate}), is
of infinite order. For $r=3,$ we have
$\mathcal{E}_2(3):V^2=U^3-14410137600U+662504472576000,$ which has
rank 2. Then there are infinitely many rational points on
$\mathcal{E}_r$ for $r\neq0,\pm1$.

Hence, there are infinitely many rational points on
$\mathcal{E}_{2}(r)$ and $\mathcal{C}_2$, and then for
\[f(x)=\frac{(x+b)(2rx-br^2+b)(2rx+br^2-b)}{4r^2x}\] Eq. (\ref{Eq11}) has infinitely many rational
solutions.

This completes the proof of Theorem 1.2 for Eq. (\ref{Eq11}). The
same method can be used to give a proof for Eq. (\ref{Eq12}).
\end{proof}

\section{Some related questions}

We have studied the rational parametric solutions of Eqs.
(\ref{Eq11}) and (\ref{Eq12}) for $f(x)=x+b+c/x$,~$x^2+bx+c+d/x$ and
$x+b+c/x+d/x^2$, but we don't get the same results for other Laurent
polynomials.

\begin{question}
For $f(x)=ax^3+bx^2+cx+d+e/x$ or $ax+b+c/x+d/x^2+e/x^3$, do Eqs.
(\ref{Eq11}) and (\ref{Eq12}) have rational solutions? If they have,
are there infinitely many?
\end{question}

Finding the integer solutions of Eqs. (\ref{Eq11}) and (\ref{Eq12})
is also an interesting question. Here we list some nontrivial
integer solutions of Eq. (\ref{Eq11}) for $f(x)=x+1+c/x$ with $|x| <
|y|\leq 300$ in Table 1.

\[\begin{tabular}{c|c}
\hline

 $c$ & $(x,y,z)$ \\

\hline

 $4$ & $(-4, -2, 5),(-2, -1, 5)$  \\
\hline

 $6$ & $(-10, 5, 12),(-6, 1, 10),(-6, 6, 10),(-1, 1, 10),$\\&$(-1, 6, 10),(1, 2, 10),(1, 3, 10),(2, 6, 10),(3, 6, 10)$  \\
 \hline

 $10$ & $(-5, 2, 10),(-5, 5, 10),(-2, 2, 10),(-2, 5, 10)$  \\
 \hline

 $12$ & $(-12, 2, 15),(-12, 6, 15),(-4, 3, 10),(-4, 4, 10),$\\&$(-3, 3, 10),(-3, 4, 10),(-1, 2, 15),(-1, 6, 15)$  \\

 \hline

 $18$ & $(5, 10, 16)$  \\

\hline

 $22$ & $(-110, -55, 122),(-10, -5, 14)$  \\
 \hline
\end{tabular}
\]
\begin{center}Table~~1. Some integer solutions of Eq. (\ref{Eq11})
for $f(x)=x+1+c/x$ with~$|x| < |y|\leq 300$ \end{center}

By some calculations, we see that a lot of Laurent polynomials
$f(x)=ax+b+c/x$ such that Eqs. (\ref{Eq11}) and (\ref{Eq12}) have
integer solutions, but we can't solve the following question

\begin{question} Does there exist a Laurent polynomial $f(x)=ax+b+c/x$ such that Eqs. (\ref{Eq11}) and
(\ref{Eq12}) have infinitely many integer solutions?
\end{question}

In a recent paper of Sz. Tengely and M. Ulas \cite{Tengely-Ulas},
they studied the integer solutions of the Diophantine equations
\begin{equation}\label{Eq31}
z^2=f(x)^2\pm g(y)^2
\end{equation}
for different $f$ and $g$. Motivated by this, we have

\begin{question}
For $f(x)=ax^2+bx+c+d/x$ and $g(x)=ex^2+fx+g+h/x$ (or other simple
Laurent polynomials), does Eq. (\ref{Eq31}) have rational solutions?
If they have, are there infinitely many?
\end{question}

\end{document}